\documentclass[reqno,a4paper]{amsart}

\usepackage{amssymb}
\usepackage{multicol}
\usepackage{amsfonts}
\usepackage{mathtools}
\usepackage{graphicx}
\usepackage{amsmath, epsfig}

\newtheorem{theorem}{Theorem}

\newtheorem{lemma}{Lemma}
\newtheorem{corollary}{Corollary}

\newtheorem{example}{Example}
\newtheorem{proposition}{Proposition}
\newcommand{\sym}{\mathcal{S}}

\newcommand{\A}{\mathcal{A}}
\newcommand{\B}{\mathcal{B}}
\newcommand{\Bclass}{\B}
\def\bangle{\atopwithdelims \langle \rangle}
\newcommand{\euler}[2]{\ensuremath{{{#1} \bangle {#2}}}}

\newcommand{\BB}{\mathbf B}
\newcommand{\NN}{\mathbb{N}}
\newcommand{\ZZ}{\mathbb{Z}}
\newcommand{\Top}{\mathsf{Top}}
\newcommand{\symfun}{\mathsf{S}}
\newcommand{\Euler}{A}

\DeclareMathOperator{\ssort}{\mathsf{stack}}
\DeclareMathOperator{\bsort}{\mathsf{bubble}}
\DeclareMathOperator{\bsc}{bsc}
\DeclareMathOperator{\st}{st}
\DeclareMathOperator{\Des}{Des}
\DeclareMathOperator{\des}{des}

\DeclareMathOperator{\maxdrop}{maxdrop}
\DeclareMathOperator{\tl}{\mathit{t}}

\newcommand{\ds}{\displaystyle}

\begin{document}
\title[Descent polynomials for permutations with bounded drop]
      {Descent polynomials for permutations with bounded drop size}

\author[F. Chung]{Fan Chung}
\author[A. Claesson]{Anders Claesson}
\author[M. Dukes]{Mark Dukes}
\author[R. Graham]{Ronald Graham}

\thanks{AC and MD were supported by grant no. 090038011 from the
  Icelandic Research Fund.}

\address{F. Chung and R. Graham:
University of California at San Diego, La Jolla CA 92093, USA}
\address{A. Claesson: The Mathematics Institute,
  School of Computer Science, 
  Reykjavik University, 103 Reykjavik, Iceland} 
\address{M. Dukes: Science Institute, University of Iceland, 107
Reykjavik, Iceland}

\begin{abstract}
  Motivated by juggling sequences and bubble sort, we examine
  permutations on the set $\{1, 2, \dots, n\}$ with $d$ descents and
  maximum drop size $k$. We give explicit formulas for enumerating
  such permutations for given integers $k$ and $d$.  We also derive
  the related generating functions and prove unimodality and symmetry
  of the coefficients.
\end{abstract}

\maketitle
\thispagestyle{empty}


\section{Introduction}
There have been extensive studies of various statistics on $\sym_n$,
the set of all permutations of $\{ 1,2, \dots, n\}$.  For a
permutation $\pi$ in $\sym_n$, we say that $\pi$ has a \emph{drop} at
$i$ if $\pi_i < i$ and that the \emph{drop size} is $i-\pi_i$.  We say
that $\pi$ has a \emph{descent} at $i$ if $\pi_i > \pi_{i+1}$. One of
the earliest results \cite{mac} in permutation statistics states that
the number of permutations in $\sym_n$ with $k$ drops equals the
number of permutations with $k$ descents. A concept closely related to
drops is that of \emph{excedances}, which is just a drop of the
inverse permutation. In this paper we focus on drops instead of
excedances because of their connection with our motivating
applications concerning bubble sort and juggling sequences.

Other statistics on a permutation $\pi$ include such things as the
number of \emph{inversions}, that is, $|\{(i,j) : i < j,\; \pi_i >
\pi_j\}|$, and the \emph{major index} of $\pi$ (i.e., the sum of $i$
for which a descent occurs). The enumeration of and generating
functions for these statistics can be traced back to the work of
Rodrigues in 1839 \cite{rodrigues} but was mainly influenced by
McMahon's treatise in 1915 \cite{mac}. There is an extensive
literature studying the distribution of the above statistics and their
$q$-analogs (see for example Foata and Han \cite{foata_han} or the
papers of Shareshian and Wachs \cite{wachs,wachs2} for more recent
developments). 
As noted above, the drop statistic that we study is closely related to the 
excedances statistic. The distribution of the bivariate 
statistics $(\mbox{descents},\mbox{excedances})$ can be found in 
Foata and Han~\cite[equations (1.15) and (1.16)]{foata_han_fix}.

This joint work originated from its connection with a paper \cite{jug}
on sequences that can be translated into juggling patterns. The set of
juggling sequences of period $n$ containing a specific state, called
the ground state, corresponds to the set $\B_{n,k}$ of permutations in
$\sym_n$ with drops of size at most $k$. As it turns out, $\B_{n,k}$
can also be associated with the set of permutations that can be sorted
by $k$ operations of bubble sort. These connections will be further
described in the next section. We note that the {\it maxdrop}
statistic has not been treated in the literature as extensively as
many other statistics in permutations. As far as we know, this is the
first time that the distribution of descents with respect to maxdrop
has been determined.

First we give some definitions concerning the statistics and
polynomials that we examine. Given a permutation $\pi$ in $\sym_n$,
let $\Des(\pi)$ denote the descent set, $\{1\leq i<n: \pi_i >
\pi_{i+1}\}$, of $\pi$ and let $\des(\pi)=|\Des(\pi)|$ be the number
of descents. We use $\maxdrop(\pi)$ to denote the value of the maximum
drop (or maxdrop) of $\pi$,
\begin{equation*}
  \maxdrop (\pi) = \max\{\,i - \pi(i) : 1\leq i\leq n\,\}.
\end{equation*}
Let $\B_{n,k}=\{\pi \in \sym_n : \maxdrop(\pi) \leq k\}$. It is
known, and also easy to show, that $|\B_{n,k}| = k!(k+1)^{n-k}$; e.g., see
\cite[Thm. 1]{jug} or \cite[p. 108]{knuth}. Let
\begin{equation*}
  b_{n,k}(r) = |\{\pi \in \B_{n,k} : \des(\pi) = r \}|,
\end{equation*}
and define the ($k$-maxdrop-restricted) descent polynomial
\begin{equation*}
  B_{n,k} (x) = \sum_{r \geq 0} b_{n,k}(r) x^r = \sum_{\pi \in \B_{n,k}}x^{\des(\pi)}.  
\end{equation*}

Examining the case of $k=2$, we discovered that the coefficients $b_{n,2}(r)$
of $B_{n,2}(x)$ appear to be given by every \emph{third} coefficient of
the simple polynomial
\begin{equation*}
  (1+x^2)(1+x+x^2)^{n-1}.
\end{equation*}
Looking at the next two cases, $k=3$ and $k=4$, yielded more
mysterious polynomials: $b_{n,3}(r)$ appeared to be every fourth
coefficient of
\begin{equation*}
  (1+x^2+2x^3+x^4+x^6)(1+x+x^2+x^3)^{n-2}
\end{equation*}
and $b_{n,4}(r)$ every fifth coefficient of 
\begin{equation*}
  (1+x^2+2x^3+4x^4+4x^5+4x^7+4x^8+2x^9+x^{10}+x^{12})(1+x+x^2+x^3+x^4)^{n-3}.
\end{equation*}
After a fierce battle with these polynomials, we were able to show that
$b_{n,k}(r)$ is the coefficient of $ u^{r(k+1)}$ in the polynomial
\begin{equation}\label{closed_1}
  P_k(u) \left(1+u+\dots + u^k \right)^{n-k} 
\end{equation}
where
\begin{equation}\label{closed_2}
  P_k(u) 
  = \sum_{j=0}^k \Euler_{k-j}(u^{k+1})(u^{k+1}-1)^{j}\sum_{i=j}^k\binom{i}{j}u^{-i},
\end{equation}
and $\Euler_k$ denotes the $k$th Eulerian polynomial (defined in the
next section). Further to this, we give an expression for the
generating function $\BB_k(z,y) = \sum_{n \geq 0}B_{n,k}(y)z^n$,
namely
\begin{equation*}
  \BB_k(z,y) = \dfrac{\ds{ 1+\sum_{t=1}^k \left(
    \Euler_t(y) - \sum_{i=1}^t \binom{k+1}{i} (y-1)^{i-1} \Euler_{t-i}(y) 
    \right)z^t }}{\ds{ 1 - \sum_{i=1}^{k+1}\binom{k+1}{i}z^i (y-1)^{i-1} }}.
\end{equation*}

We also give some alternative formulations for $P_k$ which lead to
some identities involving Eulerian numbers as well as proving the
symmetry and unimodality of the polynomials $B_{n,k}(x)$.
  
Many questions remain. For example, is there a more natural bijective
proof for the formulas that we have derived for $B_{n,k}$ and $\BB_k$?
Why do permutations that are $k$-bubble sortable define the
aforementioned juggling sequences?
 
\section{Descent polynomials, bubble sort and juggling sequences}
We first state some standard notation. The polynomial
\begin{equation*}
  \Euler_n(x) = \sum_{\pi\in\sym_n} x^{\des(\pi)}
\end{equation*}
is called the $n$th \emph{Eulerian polynomial}. For instance,
$\Euler_0(x)=\Euler_1(x)=1$ and $\Euler_2(x)=1+x$. Note that
$B_{n,k}(x) = \Euler_n(x)$ for $k \geq n-1$, since $\maxdrop(\pi) \leq
n-1$ for all $\pi \in \sym_n$.  The coefficient of $x^k$ in
$\Euler_n(x)$ is denoted $\euler{n}{k}$ and is called an \emph{Eulerian
number}. It is well known that (\cite{concrete})
\begin{equation}\label{euler_gf}
  \frac{1-w}{e^{(w-1)z }- w }
  = \sum_{k,n \geq 0} \euler{n}{k} w^k \frac{z^n}{n!}.
\end{equation}
The Eulerian numbers are also known to be given explicitly as
(\cite{euler, concrete})
\begin{equation*}
  \euler{n}{k} = \sum_{i=0}^n \binom{n+1}{i}(k+1-i)^n (-1)^i .
\end{equation*}
 
We define the operator $\bsort$ which acts recursively on permutations
via 
\begin{equation*}
  \bsort(LnR)=\bsort(L)Rn.  
\end{equation*}
In other words, to apply $\bsort$ to a permutation $\pi$ in $\sym_n$,
we split $\pi$ into (possibly empty) blocks $L$ and $R$ to the left
and right, respectively, of the largest element of $\pi$ (which
initially is $n$), interchange $n$ and $R$, and then recursively apply
this procedure to $L$. We will use the convention that
$\bsort(\emptyset) = \emptyset$; here $\emptyset$ denotes the
empty permutation. This operator corresponds to one \emph{pass} of the
classical bubble sort operation. Several interesting results on the
analysis of bubble sort can be found in
Knuth~\cite[pp. 106--110]{knuth}. We define the \emph{bubble sort
  complexity} of $\pi$ as
\begin{equation*}
  \bsc(\pi) = \min \{ k: \bsort^k(\pi)=\mbox{id}\},
\end{equation*}
the number of times $\bsort$ must be applied to $\pi$ to give the
identity permutation. The following lemma is easy to prove using
induction.

\begin{lemma}\label{bubblesort}
  $\mathrm{(i)}$ For all permutations $\pi$ we have $\maxdrop(\pi) =
  \bsc(\pi)$.\\ 
  $\mathrm{(ii)}$ The bubble sort
  operator maps $\B_{n,k}$ to $\B_{n,k-1}$.
\end{lemma}

The analysis of algorithms similar to bubble sort has been
instrumental in generating interesting research.  For example, the
analysis of stack sort in Knuth\cite[pp. 242--243]{knu} gave rise to
the area of pattern avoiding permutations.  The stack sort operator
$\ssort$ is defined by $\ssort(LnR)=\ssort(L)\ssort(R)n$. We see below
that stack sort is at least as efficient as bubble sort.

\begin{lemma}
\label{stack}
For all $\pi \in \sym_n$, if $\bsort^k(\pi) = \mbox{id}$ then
$\ssort^k(\pi)=\mbox{id}$.
\end{lemma}

\begin{proof} 
  The proof of Lemma \ref{stack} follows from the following claim:
  \smallskip

  {\it If $A=a_1a_2 \dots a_n=LmR$ is a sequence of
    distinct positive integers and $m=\max_i a_i$, then
    either $\maxdrop(A) = 1-a_1$ or $\maxdrop(\ssort(A)) \leq
    \maxdrop(A) -1$.}\smallskip

  The Claim is certainly true for $n=1$. Suppose the claim is true for
  $n' < n$.  If $\maxdrop(A)=1-a_1$, we are done. We may assume that
  $\maxdrop(A) > 1-a_1$.  This implies that the maxdrop of $A$ does
  not occur at the entry where $m$ is located.  For the $i$th entry in
  $\ssort(L)$, the maxdrop of $\ssort(L)$ at $i$ is reduced by one by
  induction.  For the $j$th entry in $R$, the maxdrop of the
  corresponding entry in $\ssort(A)$ is reduced by $1$.  Thus, the
  claim is proved by induction.
\end{proof}

The class of permutations $\B_{n,k}$ appears in a recent paper
\cite{jug} on enumerating juggling patterns that are usually 
called \emph{siteswaps} by (mathematically inclined) jugglers. 
Suppose a juggler throws a ball at time $i$ so that the ball will be
in the air for a time $t_i$ before landing at time $t_i +i$. Instead
of an infinite sequence, we will consider periodic patterns, denoted
by $T=(t_1, t_2, \dots, t_n)$.  A \emph{juggling sequence} is just one
in which two balls never land at the same time. It is not hard to show
\cite{jug0} that a necessary and sufficient condition for a sequence
to be a juggling sequence is that all the values $t_i+i \pmod n$ are
distinct.  In particular, it follows that that the average of $t_i$ is
just the numbers of balls being juggled. Here is an example: 

If $T=(3,5,0,2,0)$ then at time 1 a ball is thrown that will land at
time $1+3=4$. At time 2 a ball is thrown that will land at time
$2+5=7$. At time 3 a ball is thrown that will land at time
$3+0=3$. Alternatively one can say that no ball is thrown at time 3.
This is represented in the following diagram. \ \\[0.6em]
\centerline{\scalebox{0.75}{\includegraphics{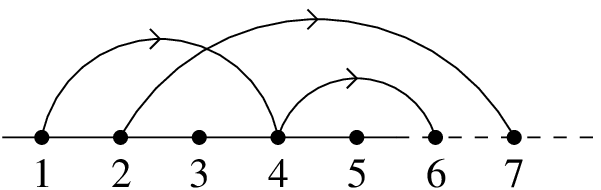}}}
Repeating this for all intervals of length 5 gives \ \\[0.9em]
\centerline{\scalebox{0.75}{\includegraphics{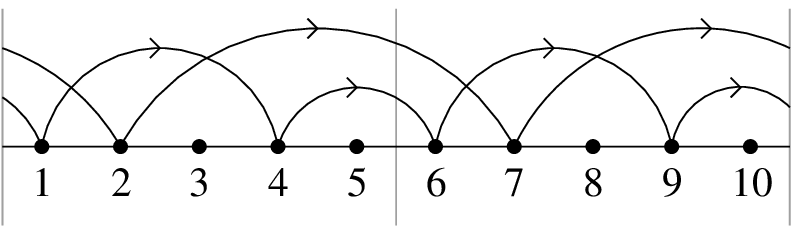}}}
 
For a given juggling sequence, it is often possible to further
decompose into shorter juggling sequences, called \emph{primitive
  juggling sequences}, which themselves cannot be further decomposed.
These primitive juggling sequences act as basic building blocks for
juggling sequences \cite{jug}.  However, in the other direction, it is
not always possible to combine primitive juggling sequences into a
longer juggling sequence. Nevertheless, if primitive juggling
sequences share a common \emph{state} (which one can think of as a
\emph{landing schedule}), we then can combine them to form a longer
and more complicated juggling sequences.  In \cite{jug}, primitive
juggling sequences associated with a specified state are enumerated.
Here we mention the related fact concerning $\B_{n,k}$:\smallskip

{\it There is a bijection mapping permutations in $\B_{n,k}$ to
  primitive juggling sequences of period $n$ with $k$ balls that all
  share a certain state, called the ground state.}\smallskip

The bijection maps $\pi$ to $\phi(\pi)= (t_1, \dots, t_n) $ with $t_i
= k-i+\pi_i$. As a consequence of the above fact and Lemma
\ref{bubblesort}, we can use bubble sort to transform a juggling
sequence using $k$ balls to a juggling sequence using $k-1$ balls.

To make this more precise, let $T=(t_1,\dots, t_n)$ be a juggling
sequence that corresponds to $\pi \in \Bclass_{n,k}$, and suppose that
$T'=(s_1,\dots , s_n)$ is the juggling sequence that corresponds to
$\bsort(\pi)$.  Assume that the ball $B$ thrown at time $j$ is the one
that lands latest out of all the $n$ throws.  In other words, $t_j+j$
is the largest element in $\{t_i+i\}_{i=1}^n$. Now, write $T=L t_j R $
where $L=(t_1,\dots,t_{j-1})$ and $R=(t_{j+1},\dots,t_n)$. Then we
have
\begin{equation*}
  T' = f_k(T) = f_k(L)R\hspace{0.7pt}s,  
\end{equation*}
where $s = t_j + j - (n+1)$. In other words, we have removed the ball
$B$ thrown at time $j$ and thus throw all balls after time $j$ one
time unit sooner.  Then at time $n$ we throw the ball B so that it
lands one time unit sooner than it would have originally landed.  Then
we repeat this procedure to all the balls thrown before time $j$.


\section{The polynomials $B_{n,k}(y)$}
In this section we will characterise the polynomials $B_{n,k}(y)$.
This is done by first finding a recurrence for the polynomials and
then solving the recurrence by exploiting some aspects of their
associated characteristic polynomials. The latter step is quite
involved and so we present the special case dealing with $B_{n,4}(y)$
first.

\subsection{Deriving the recurrence for $B_{n,k}$}
We will derive  the following recurrence for $B_{n,k}(y)$.
\begin{theorem}\label{b_rec_thm}
  For $n \geq 0$,
  \begin{equation}\label{rec}
    B_{n+k+1,k} (y) 
    = \sum_{i=1}^{k+1} \binom{k+1}{i} (y-1)^{i-1} B_{n+k+1-i,k} (y)
  \end{equation}
  with the initial conditions
  \begin{equation*}
    B_{i,k}(y)  = \Euler_i(y),\quad 0 \leq i \leq k.
  \end{equation*}
\end{theorem}

We use the notation $[a,b] =\{i\in\ZZ :a\leq i\leq b\}$ and
$[b]=[1,b]$. Let $A=\{a_1,\dots,a_n\}$ with $a_1<\dots<a_n$ be any
finite subset of $\NN$. The \emph{standardization} of a permutation
$\pi$ on $A$ is the permutation $\st(\pi)$ on $[n]$ obtained from
$\pi$ by replacing the integer $a_i$ with the integer $i$. Thus $\pi$
and $\st(\pi)$ are order isomorphic. For example, $\st(19452) =
15342$. If the set $A$ is fixed, the inverse of the standardization
map is well defined, and we denote it by $\st^{-1}_A(\sigma)$; for
instance, with $A=\{1,2,4,5,9\}$, we have $\st^{-1}_A(15342)=19452$.
Note that $\st$ and $\st^{-1}_A$ each preserve the descent set.

For any set $S \subseteq [n-1]$ we define
$\A_{n,k}(S) = \{ \pi\in\Bclass_{n,k} : \Des(\pi) \supseteq S \}$
and
\begin{equation*}
  t_n(S) = \max\{ i\in\NN : [n-i,n-1] \subseteq S\}.
\end{equation*}
Note that $\tl_n(S)=0$ in the case that $n-1$ is not a member of
$S$. Now, for any permutation $\pi=\pi_1\dots\pi_n$ in $\A_{n,k}(S)$ define
\begin{equation*}
  f(\pi) = (\sigma, X),\;\text{ where }
  \sigma=\st(\pi_1\dots\pi_{n-i-1}),
  X=\{\pi_{n-i},\dots ,\pi_n\} \text{ and } i=\tl_n(S).  
\end{equation*}

\begin{example} 
  Let $S=\{3,7,8\}$, and choose the permutation $\pi = 138425976$ in
  $\A_{9,3}(S)$. Notice that $\Des(\pi)=\{3,4,7,8\} \supset S$.  Now
  $\tl_9(S)=2$. This gives $f(\pi) = (\sigma, X)$ where $\sigma =
  \st(138425) = 136425$ and $X=\{\pi_7, \pi_8,\pi_9\}=\{6,7,9\}$.
  Hence $f(138425976)=(136425, \{6,7,9\})$.
\end{example}

\begin{lemma}
  For any $\pi$ in $\A_{n,k}(S)$, the image $f(\pi)$ is in the
  Cartesian product
  \begin{equation*}
    \A_{n-i-1,k}(S\cap [n-\tl_n(S)-2])\times\binom{[n-k,n]}{\tl_n(S)+1},
  \end{equation*}
  where $\binom{X}{m}$ denotes that set of all $m$-element subsets
  of the set $X$.
\end{lemma}

\begin{proof}
  Given $\pi\in\A_{n,k}(S)$, let $f(\pi)=(\sigma,X)$. Suppose
  $i=\tl_n(S)$. Then there are descents at positions $n-i, \dots ,
  n-1$ (this is an empty sequence in case $i=0$). Thus
  \begin{equation*}
    n\geq\pi_{n-i}>\pi_{n-i+1}>\dots > \pi_{n-1}>\pi_n\geq n-k,  
  \end{equation*}
  where the last inequality follows from the assumption that
  $\maxdrop(\pi)\leq k$. Hence $X$ is an $(i+1)$-element subset of
  $[n-k,n]$, as claimed. Clearly $\sigma\in\sym_{n-i-1}$.

  Next we shall show that $\sigma$ is in $\A_{n-i-1,k}$.  Notice that
  the entries of $(\pi_1,\dots, \pi_{n-i-1})$ that do not change
  under standardization are those $\pi_{\ell}$ which are
  $<\pi_n$. Since these values remain unchanged, the values
  $\ell-\pi_{\ell}$ are also unchanged and are thus $\leq k$.
  
  Let $(\pi_{a(1)},\dots, \pi_{a(m)})$ be the subsequence of values
  which are $>\pi_n$. The smallest value that any of these may take
  after standardization is $\pi_n \geq n-k$.  So $\sigma_{a(j)} \geq
  \pi_n \geq n-k$ for all $j \in [1,m]$.  Thus $a(j)-\sigma_{a(j)}
  \leq a(j) - (n-k) = k-(n-a(j)) \leq k$ for all $j \in [1,m]$.
  Therefore $\ell - \sigma_{\ell} \leq k$ for all $\ell \in [1,n-i-1]$
  and so $\sigma \in \A_{n-i-1,k}$.

  The descent set is preserved under standardization, and consequently
  $\sigma$ is in $\A_{n-i-1,k}(S\cap [n-i-2])$, as claimed.
\end{proof}

We now define a function $g$ which will be shown to be the inverse of
$f$. Let $\pi$ be a permutation in $\A_{m,k}(T)$, where $T$ is a
subset of $[m-1]$. We will add $i+1$ elements to $\pi$ to yield a new
permutation $\sigma$ in $\A_{m+i,k}(T \cup [m+1,m+i])$.  Choose any
$(i+1)$-element subset $X$ of the interval $[m+i+1-k,m+i+1]$, and let
us write $X=\{x_1,\dots , x_{i+1}\}$, where $x_1\leq\dots\leq
x_{i+1}$. Define
\begin{equation*}
  g(\pi,X)=\st^{-1}_V (\pi_1\dots\pi_m) \, x_{i+1}x_i\dots x_1,\,
  \text{ where }V=[m+i+1]\setminus X.
\end{equation*}

\begin{example}
  Let $T=\{1\}$, and choose the permutation $\pi=3142$ in $\A_{4,3}
  (T)$. Notice that $\Des(\pi)=\{1,3\}\supseteq T$. Choose $i=2$ and
  select a subset $X$ from $[4+2+1-3,4+2+1]=\{4,5,6,7\}$ of size
  $i+1=3$. Let us select $X=\{4,6,7\}$. Now we have
  $g(\pi,X) =\st^{-1}_V(3142)\,764=3152764$, 
  where $V$ is the set $[4+2+1]\setminus \{4,6,7\} = \{1,2,3,5\}$.
\end{example}

\begin{lemma}
  If $(\pi,X)$ is in the Cartesian product 
  \begin{equation*}
    \A_{m,k}(T)\times\binom{[m+i+1-k,m+i+1]}{i+1}
  \end{equation*}
  for some $i>0$ then $g(\pi,X)$ is in
  \begin{equation*}
    \A_{m+i+1,k}(T\cup [m+1, m+i]).
  \end{equation*}
\end{lemma}

\begin{proof}
  Let $\sigma=g(\pi,X)$. For the first $m$ elements of $\sigma$, since
  $\sigma_j \geq \pi_j$ for all $1\leq j \leq m$, we have $j-\sigma_j
  \leq j-\pi_j$
  which gives 
  \begin{equation*}
    \max\{j-\sigma_j: j\in[m]\}\leq\max\{j-\pi_j : j\in[m]\}\leq k.
  \end{equation*}
  The final $i+1$ elements of $\sigma$ are decreasing so the
  $\maxdrop$ of these elements will be the $\maxdrop$ of the final
  element, 
  \begin{equation*}
    m+i+1-\sigma_{m+i+1} = m+i+1 - x_1 \leq m+i+1-(m+i+1-k) = k.
  \end{equation*}
  Thus $\maxdrop(\sigma) \leq k$ and so $\sigma\in\Bclass_{m+i+1,k}$.
  The descents of $\sigma$ will be in the set $T \cup
  [m+1,m+i]$ since descents are preserved under standardization 
  and the final $i+1$
  elements of $\sigma$ are listed in decreasing order. Hence 
  $\sigma\in A_{m+i+1,k}(T\cup [m+1, m+i])$, as claimed.
\end{proof}

\begin{lemma}
  The function $f$ is a bijection, and $g$ is its inverse.
\end{lemma}

\begin{proof}
  Given any $(\sigma,X) \in \binom{[m+j+1-k,m+j+1]}{j+1} 
  \times\A_{m,k}(T)$ where $T \subseteq [m-1]$, let
  $\pi=g(\sigma,X)$. We have
  \begin{equation*}
    \pi=\st^{-1}_{[m+j+1]\setminus X}(\sigma_1\dots\sigma_m)
    \,x_{j+1}x_j\dots x_1 \in \A_{m+j+1,k}(T\cup[m+1,m+j]),    
  \end{equation*}
  where $X=\{x_1,\dots, x_{j+1}\}$ and $x_1\leq\dots\leq x_{j+1}$.
  Let $S=T\cup [m+1,m+j]$. Clearly $\Des(g(\sigma,X)) \supseteq S$
  and $i=\tl_{m+j+1}(S)=j$. So $f(g(\sigma,X)) = (\tau,Y)$ where
  \begin{equation*}
    Y=\{\pi_{m+1},\dots , \pi_{m+j+1}\} = \{x_1,\dots, x_{j+1}\}= X
  \end{equation*}
  and
  \begin{equation*}
    \tau = 
    \st(\st^{-1}_{[1,m+j+1]\setminus X}(\sigma_1\dots\sigma_m)) = \sigma.
  \end{equation*}
  Hence
  $f(g(\sigma,X))=(\sigma,X)$.
  Given $\pi\in\A_{n,k}(S)$, let $f(\pi)=(\sigma,X)$ with
  $X=\{x_1,\dots , x_{i+1}\}$ and 
  $\sigma = \st(\pi_1\dots\pi_{n-(i+1)})$.
  We have
  \begin{align*}
    g(\sigma,X) 
    &=\st^{-1}_{[1,n]\setminus X} (\sigma_1\dots\sigma_{n-i-1})
    x_{i+1}\dots x_1\\
    &= \st^{-1}_{[1,n]\setminus X}(\st(\pi_1\dots\pi_{n-i-1}))
    \pi_{n-i}\dots \pi_n \\
    &= \pi_1\dots \pi_{n-i-1}\pi_{n-i}\dots \pi_n \\
    &= \pi.
  \end{align*}
  Hence $g(f(\pi))=\pi$.
\end{proof}

\begin{corollary}\label{rec_corol}
  Let $a_{n,k}(S) = |\A_{n,k}(S)|$ and $i=\tl_n(S)$. Then
  \begin{equation*}
    a_{n,k}(S) = \binom{k+1}{i+1} a_{n-(i+1),k} (S\cap [1,n-(i+1)]).
  \end{equation*}
\end{corollary}

\begin{proposition}\label{prop_rec}
  For all $n\geq0$,
  $$\Bclass_{n,k}(y+1) = 
  \sum_{i=1}^{k+1} {k+1 \choose i} y^{i-1} \Bclass_{n-i,k} (y+1).
  $$
\end{proposition}

\begin{proof}
  Notice that
  \begin{align*}
    \Bclass_{n,k}(y+1) 
    &= \sum_{\pi \in \Bclass_{n,k}} (y+1)^{\des(\pi)} \\
    &= \sum_{\pi \in \Bclass_{n,k}} \sum_{i=0}^{\des(\pi)} \binom{\des(\pi)}{i} y^i \\
    &= \sum_{\pi \in \Bclass_{n,k}} \sum_{S\subseteq \Des(\pi)} y^{|S|}\\
    &= \sum_{S\subseteq [n-1]} y^{|S|} \sum_{\pi \in \A_{n,k}(S)}  1
    = \sum_{S\subseteq [n-1]} y^{|S|} a_{n,k}(S).
  \end{align*}
  From Corollary \ref{rec_corol}, multiply both sides by $y^{|S|}$ and sum
  over all $S \subseteq [n-1]$. We have
  \begin{align*}
    \Bclass_{n,k}(y+1) 
    &= \sum_{S \subseteq [n-1]} y^{|S|} \binom{k+1}{\tl_n(S)+1} a_{n-(\tl_n(S)+1),k} (S \cap [n-(\tl_n(S)+2)]) \\
    &= \sum_{i \geq 0 } \sum_{S \subseteq [n-1]\atop \tl_n(S)=i} y^i y^{|S|-i} \binom{k+1}{i+1} a_{n-(i+1),k} (S \cap [n-(i+2)]) \\
    &= \sum_{i \geq 0 } \binom{k+1}{i+1} y^i \sum_{S \subseteq [n-1]\atop \tl_n(S)=i}  a_{n-(i+1),k} (S \cap [n-(i+2)]) y^{|S|-i}\\
    &= \sum_{i \geq 0 } \binom{k+1}{i+1} y^i \sum_{S \subseteq [n-(i+1)]}  a_{n-(i+1),k} (S) y^{|S|}\\
    &= \sum_{i \geq 0 } \binom{k+1}{i+1} y^i \Bclass_{n-(i+1),k} (y+1) \\
    &= \sum_{i \geq 1 } \binom{k+1}{i} y^{i-1} \Bclass_{n-i,k} (y+1).
  \end{align*}
\end{proof}

\begin{proof}[Proof of Theorem \ref{b_rec_thm}]
  Replacing $n$ and $y$ by $n+k+1$ and $y-1$, respectively, in
  Proposition \ref{prop_rec} yields the recurrence (\ref{rec}):
  \begin{equation*}
    B_{n+k+1,k}(y) = \sum_{i=1}^{k+1} \binom{k+1}{i} (y-1)^{i-1} B_{n+k+1-i,k} (y)  
  \end{equation*}
  for $n \geq 0$, with the initial conditions
  $B_{i,k}(y)  = \Euler_i(y),\,0 \leq i \leq k$.
\end{proof}
Consequently, by multiplying the above recurrence by $z^n$ and summing
over all $n\geq 0$, we have the generating function $\BB_k(z,y)$:
\begin{equation}\label{BB_gen}
  \BB_k(z,y) = \dfrac{\ds{ 1+\sum_{t=1}^k \left(
    \Euler_t(y) - \sum_{i=1}^t \binom{k+1}{i} (y-1)^{i-1} \Euler_{t-i}(y) 
    \right)z^t }}{\ds{ 1 - \sum_{i=1}^{k+1}\binom{k+1}{i}z^i (y-1)^{i-1} }}.
\end{equation}

\subsection{Solving the recurrence for $B_{n,4}$.}
Before we proceed to solve the recurrence for $B_{n,k}$, we first
examine the special case of $k=4$ which is quite illuminating. We note
that the characteristic polynomial for the recurrence for $B_{n,4}$ is
\begin{align*}
  h(z) 
  &= z^5 - 5z^4+10(1-y)z^3 - 10(1-y)^2z^2 + 5(1-y)^3z -(1-y)^4\\
  &= \frac{(z-1+y)^5 -yz^5}{1-y}.
\end{align*}
Substituting $y = t^5$ in the expression above, we see that the
roots of $h(z)$ are just 
\begin{equation*}
  \rho_j(t) = \frac{1-t^5}{1-\omega^j t},\quad 0\leq j \leq 4,
\end{equation*}
where $\omega = \exp(\frac{2 \pi \mathfrak{i}}{5})$ is
a primitive $5th$ root of unity. Hence, the general term for
$B_{n,4}(t)$ can written as
\begin{equation*}
  B_{n,4}(t) = \sum_{i=0}^4 \alpha_i(t) \rho_i^n(t) 
\end{equation*}
where the $\alpha_i(t)$ are appropriately chosen coefficients
(polynomials in $t$). To determine the $\alpha_i(t)$ we need to solve
the following system of linear equations:
\begin{equation*}
  \sum_{i=0}^4 \alpha_i(t) \rho_i^j(t) 
  = B_{j,4}(t) = \Euler_j(t^5), \, 0 \leq j \leq 4.  
\end{equation*}
Thus, $\alpha_i(t)$ can be expressed as the ratio
$N_{4,i+1}(t)/D_4(t)$ of two determinants. The denominator
$D_4(t)$ is just a standard Vandermonde determinant whose $(i+1,j+1)$
entry is $\rho_i^j(t)$. The numerator $N_{4,i+1}(t)$ is formed from
$D_4(t)$ by replacing the elements $\rho_i^j(t)$ in the $(i+1)$st row
by $\Euler_j(t^5)$.  A quick computation (using the symbolic
computation package Maple) gives:
\begin{align*}
 D_4(t)    &= 25 \sqrt 5 \,(1-t^5)^6 t^{10}; \\
 N_{4,1}(t) &=5 \sqrt 5 \,(t^{12}+t^{10}+2t^9+4t^8+4t^7+4t^5+4t^4+2t^3+t^2+1)
 (1-t^5)^3 (1-t)^3t^{10}
\end{align*}
and, in general, $N_{4,i+1}(t) = N_{4,1}(\omega^i t).$

Substituting the value $\alpha_0 (t) = N_{4,1}(t)/D_4(t)$ into the
first term in the expansion of $B_{n,4}$, we get
\begin{multline*}
  \alpha_0 (t)(1+t+t^2+t^3+t^4)^n \\
  =\tfrac{1}{5}(t^{12}+t^{10}+2t^9+4t^8+4t^7+4t^5+4t^4+2t^3+t^2+1)
  (1+t+t^2+t^3+t^4)^{n-3}.
\end{multline*}
Now, since the other four terms $\alpha_i(t)(1+t+t^2+t^3+t^4)^n$ arise
by replacing $t$ by $\omega ^i t$ then in the sum of all five terms,
the only powers of $t$ that survive are those which have powers which
are multiples of $5$.  Thus, we can conclude that if we write
\begin{equation*}
  (t^{12}+t^{10}+2t^9+4t^8+4t^7+4t^5+4t^4+2t^3+t^2+1)(1+t+t^2+t^3+t^4)^{n-3} 
  = \sum_r \beta(r) t^r   
\end{equation*}
then $b_{n,4}(d) = \beta(5d)$. In other words, the number of
permutations $\pi \in \B_{n,4}$ with $d$ descents is given by the
coefficient of $t^{5d}$ in the expansion of the above polynomial.
Incidentally, the corresponding results for the earlier $\B_{n,i}$ are
as follows: $b_{n,1}(d) = \beta(2d)$ in the expansion of
\begin{equation*}
  (1+t)^n = \sum_r \beta(r) t^r ,
\end{equation*}
so $b_{n,1}(d)=\binom{n}{2d}$; $b_{n,2}(d) = \beta(3d)$ in
the expansion of
\begin{equation*}
  (1+t^2)(1+t)^{n-1} = \sum_r \beta(r) t^r;
\end{equation*}
and $b_{n,3}(d) = \beta(4d)$ in the expansion of
\begin{equation*}
  (1+t^2+2t^3+t^4+t^6)(1+t)^{n-2} = \sum_r \beta(r) t^r.  
\end{equation*}

The preceding arguments have now set the stage for dealing with the
general case of $B_{n,k}$. Of course, the arguments will be somewhat
more involved but it is hoped that treating the above special case
will be a useful guide for the reader.

\subsection{Solving  the recurrence for $B_{n,k}$}
\begin{theorem}
  We have 
  $B_{n,k}(y) = \sum_d \beta_k\big((k+1)d\big) y^{(k+1)d}$, where
  \begin{equation*}
    \sum_j \beta_k(j) u^j 
    = P_k(u) \left( \frac{1-u^{k+1}}{1-u} \right)^{n-k} 
  \end{equation*}
  and
  \begin{equation*}
    P_k(u) = 
    \sum_{j=0}^k 
    \Euler_{k-j}(u^{k+1})  (u^{k+1}-1)^{j} \sum_{i=j}^k \binom{i}{j} u^{-i}.
  \end{equation*}  
\end{theorem}

\begin{proof}
  To solve (\ref{rec}) for $B_{n,k}$, we first need to compute the
  roots of the corresponding characteristic polynomial
  \begin{equation*}
    z^{n+k+1} - \sum_{i=1}^{k+1} \binom{k+1}{i} (y-1)^{i-1} z^{n+k+1-i} 
  \end{equation*}
  which can be rewritten as
  \begin{equation*}
    \big(z - (1-y)\big)^{k+1} - yz^{k+1}.  
  \end{equation*}
  Substituting $y = u^{k+1}$, this becomes
  \begin{equation*}
    \big(z-(1-u^{k+1})\big)^{k+1} - u^{k+1}z^{k+1}.  
  \end{equation*}
  The $k+1$ roots of this polynomial are easily seen to be the
  expressions
  \begin{equation*}
    \rho_i = \rho_i(u) 
    = \frac{1-u^{k+1}}{1-\omega^i u},\quad 0 \leq i \leq k,  
  \end{equation*}  
  where $\omega = \exp(\frac{2 \pi \mathfrak{i}}{k+1})$ is a primitive
  $(k+1)$st root of unity. Thus, we can express $B_{n,k}$ in the form
  \begin{equation*}
    B_{n,k}(u^{k+1}) = \sum_{i=0}^k \alpha_{k,i}(u) \rho_i^n
  \end{equation*}
  for an appropriate choice of coefficients $\alpha_{k,i}$ (which depend
  on the initial conditions).  In fact, writing down the expressions for
  the first $k+1$ $B_{n,k}$'s, we have:
  \begin{equation*}
    B_{j,k}(u^{k+1}) 
    = \sum_{i=0}^k \alpha_{k,i}(u) \rho_i^j  
    = \Euler_j(u^{k+1}),\quad \mbox{for } 0 \leq j \leq k.
  \end{equation*}
  We can solve this as a system of $k+1$ linear equations in the $k+1$
  unknown coefficients $\alpha_{k,i}(u)$ by representing the solution
  in the usual way as a ratio of two determinants.  In particular, the
  expression for $\alpha_{k,0}$ is given by
  \begin{equation}\label{alpha_0}
    \alpha_{k,0}(u) = \frac{\det R_k(u)}{\det S_k(u)}
  \end{equation}
  where $S_k(u)$ and $R_k(u)$ are ($k+1$) by ($k+1$) matrices defined by
  \begin{equation*}
    S_k(u) = (S_k(i+1,j+1))\;\text{ with }\,S_k(i+1,j+1) = \rho_i^j
  \end{equation*}
  and
  \begin{equation*}
    R_k(u) =\big(R_k(i+1,j+1)\big) \;\text{ with }\, R_k(i+1,j+1) 
     =
    \begin{cases}
      \Euler_{j}(u^{k+1})& \text{ if } i = 0, \\
      \rho_i^j          & \text{ if } i > 0.
    \end{cases}
  \end{equation*}
  Now the bottom determinant is a standard Vandermonde determinant which
  has the value
  \begin{equation*}
    \det S_k(u) =  \prod_{0 \leq i < j \leq k} (\rho_j - \rho_i).
  \end{equation*}
  The top determinant is almost a Vandermonde determinant (except for
  the first row). Its value has the form
  \begin{equation*}
    \det R_k(u) = \Top_k(u) \cdot \prod_{0 < i < j \leq k} (\rho_j - \rho_i)
  \end{equation*}
  where $\Top_k(u)$ is a polynomial in $u$ which we will soon determine.
  Hence, in the ratio (\ref{alpha_0}), the terms which do not involve
  $\rho_0$ cancel, leaving the reduced form
  \begin{equation*}
    \alpha_{k,0}(u) = \frac{\Top_k(u)}{\prod_{j > 0} (\rho_j - \rho_0)}.
  \end{equation*}
  However, we have
  \begin{align}
    \prod_{j > 0} (\rho_0 - \rho_j) 
    &= \prod_{j>0} \left( \frac{1-u^{k+1}}{1-u} - 
    \frac{1-u^{k+1}}{1-\omega^j u} \right)\nonumber \\
    &= \left( \prod_{j>0} \frac{(1-\omega^j) u}{(1-u)(1-\omega^j u)}\right)
    (1-u^{k+1})^k \nonumber \\
    &= \frac{(1-u^{k+1})^k u^k}{(1-u)^k}\cdot
    \frac{\prod_{j>0} (1-\omega^j) }{\prod_{j>0} (1-\omega^j u)} \nonumber\\
    &= \frac{(1-u^{k+1})^k u^k}{(1-u)^k}\cdot
    \frac{k+1}{\frac{1-u^{k+1}}{1-u}} 
    = \frac{(1-u^{k+1})^{k-1} u^k}{(1-u)^{k-1}} \cdot (k+1).
    \label{bot2}
  \end{align}
  On the other hand, for the top we have by standard properties of
  Vandermonde determinants:
  \begin{equation}\label{top}
    \Top_k(u) = \sum_{j \geq 0} (-1)^{k-j} \Euler_{k-j}(u^{k+1})
    \symfun_{k,j} (\rho_1,\dots, \rho_k),
  \end{equation}
  where $\symfun_{k,j}(x_1, \dots, x_k)$ is the elementary symmetric
  function of degree $j$ in the $k$ variables $x_1, \dots, x_k$ (and
  we recall that $\Euler_t(y) = \sum_{j \geq 0} \euler{t}{j} y^t$).
  Now consider the generating function
  \begin{align*}
    X_k(z) 
    &= (z-\rho_0) (z-\rho_1) \dots (z-\rho_k) \\
    &= \sum_{t \geq 0}(-1)^t\symfun_{k+1,t}(\rho_0,\dots,\rho_k)z^{k+1-t} \\ 
    &= \prod_{i \geq 0}\left(z - \frac{1-u^{k+1}}{1-\omega^i u} \right) \\
    &= \frac{\prod_{i \geq 0} (z-(1-u^{k+1}) - \omega^i uz)}
    {\prod_{i \geq 0} (1-\omega^i u)} \\
    &= \frac{\left(z-(1-u^{k+1})\right)^{k+1}-u^{k+1} z^{k+1}}{1 - u^{k+1}}.
  \end{align*}
  What we are interested in is
  \begin{align*}
    Y_k(z) = \frac{X_k(z)}{z-\rho_0}
    &= \sum_{i \geq 0} (-1)^i\symfun_{k,i}(\rho_1,\dots,\rho_k) z^{k-i}\\
    &= \frac{\left( z-(1-u^{k+1})\right)^{k+1} 
      - u^{k+1}z^{k+1}}{(1-u^{k+1})\big(z - \frac{1-u^{k+1}}{1-u}\big)}\\
    &= \frac{1-u}{1-u^{k+1}} \sum_{i=0}^k (z-1+u^{k+1})^i (uz)^{k-i} \\
    &= \frac{1-u}{1-u^{k+1}} \sum_{i=0}^k
    \sum_{j=0}^i(-1)^j\binom{i}{j}(1 - u^{k+1})^j u^{k-i}z^{k-j}\\
    &= \frac{1-u}{1-u^{k+1}} \sum_{j=0}^k(1 - u^{k+1})^j z^{k-j}
    \sum_{i = j}^k  (-1)^j \binom{i}{j}  u^{k-i}.
  \end{align*}
  Thus, by identifying the  coefficient of $z^{k-j}$, we have
  \begin{equation*}
    \symfun_{k,j}(\rho_1, \dots, \rho_k) 
    = \frac{(1-u)(1 - u^{k+1})^{j}}{1-u^{k+1}}\sum_{i=j}^k\binom{i}{j}u^{k-i}.
  \end{equation*}
  Therefore, we have
  \begin{align*}
    \Top_k(u) 
    &= \sum_{j \geq 0}(-1)^{k-j}\Euler_{k-j}(u^{k+1})
    \symfun_{k,j}(\rho_1,\dots, \rho_k) \\
    &= \sum_{j \geq 0} (-1)^{k-j}\Euler_{k-j}(u^{k+1})
    \frac{(1-u)(1 - u^{k+1})^{j}}{1-u^{k+1}}\sum_{i=j}^k\binom{i}{j}u^{k-i}.
  \end{align*}
  As a consequence, we find by (\ref{top}) and (\ref{bot2}) that
  \begin{align}\label{ex}
    \lefteqn{\alpha_{k,0}(u) \rho_0^n} \nonumber\\
    &= \frac{\Top_k(u)}{(-1)^k \prod_{j > 0} (\rho_0 - \rho_j)} 
    \left(\frac{1-u^{k+1}}{1-u}\right)^{\!n}\\
    &= \frac{1}{(k+1)u^k}\left( \frac{1-u^{k+1}}{1-u} \right)^{\!n-k}
    \sum_{j=0}^k(-1)^{j}\Euler_{k-j}(u^{k+1})  (1-u^{k+1})^{j} 
    \sum_{i=j}^k \binom{i}{j} u^{k-i}\nonumber\\
    &= \frac{1}{(k+1)} \left( \frac{1-u^{k+1}}{1-u} \right)^{\!n-k} 
    \sum_{j=0}^k\Euler_{k-j}(u^{k+1})  (u^{k+1}-1)^{j} 
    \sum_{i=j}^k \binom{i}{j} u^{-i}.\nonumber
  \end{align}
  (We should keep in mind that this expression actually is a polynomial
  in $u$.)  To determine the other coefficients $\alpha_{k,t}(u)$, we
  make the following observations.  First, for $1 \leq t \leq k$, we can
  cyclically permute the rows of the (Vandermonde) matrix $S_k(u)$ to
  form the new matrix $S_k^{(t)}(u) = (S_k^{(t)}(i+1,j+1)) =
  (\rho_i^{j+t})$. We can then form the corresponding matrix
  $R_k^{(t)}(u)$ by replacing the top row of $S_k^{(t)}(u)$ by
  $\Euler_0(u^{k+1}), \dots, \Euler_k(u^{k+1})$.  In this way, we can
  express the coefficient $\alpha_{k,t}(u)$ as:
  \begin{equation}\label{alpha_t}
    \alpha_{k,t}(u) = \frac{\det R_k^{(t)}(u)}{\det S_k^{(t)}(u)}.
  \end{equation}
  However, observe that the resulting computations for determining
  $\alpha_{k,t}(u)$ are exactly the same as those for
  $\alpha_{k,0}(u)$ where we replace $u$ by $\omega^t u$. This is
  because $\Euler_j(u^{k+1}) = \Euler_j(\omega^t u)^{k+1}$.
  Consequently,
  \begin{equation*}
    \alpha_{k,t}(u) = \alpha_{k,0}(\omega^t u).  
  \end{equation*}
  Therefore, for each $n$, the expression $\alpha_{k,t}(u) \rho_t^n$
  as a polynomial in $u$ can be obtained from $ \alpha_{k,0}(u) \rho_0^n
  $ by replacing $u$ by $\omega^t u$.  However, this implies that in the
  sum $ \sum_{t=0}^k \alpha_{k,t}(u) \rho_t^n$, the only terms that
  survive are those powers $u^m$ of $u$ which are multiples of $k+1$,
  since for $m \not\equiv 0 \pmod {k+1}$, we have $\sum_{i=0}^k
  \omega^{mi} = 0$. One the other hand, for $u^m$ with 
  $m \equiv 0 \pmod {k+1}$, we must multiply the coefficients by $k+1$ 
  since in this case, all the powers $\omega^t$ are $1$. Consequently,
  if we write
  \begin{equation}\label{final}
    (k+1) \alpha_{k,0}(u) \rho_0^n = \sum_j \beta_k(j) u^j
  \end{equation}
  then we have
  \begin{equation*}
    B_{n,k}(y) = \sum_d \beta_k\big(\,(k+1)d\,\big)\, y^{(k+1)d}.
  \end{equation*}
  In other words, the number $b_{n,k} (r)$ of permutations in $\B_{n,k}$
  having exactly $r$ descents is equal to the $k+1$ times the
  coefficient of $u^{(k+1)r}$ in (\ref{ex}).
\end{proof}
  
  If we express (\ref{ex}) (times $k+1$) in the form
  \begin{equation}
    P_k(u) \left(1+u+\dots+u^k\right)^{n-k} 
  \end{equation}
  then the first few values of $P_k(u)$ are shown below.
  \begin{equation*}
    \begin{array}{c|l}
      k & P_k(u)  \\ \hline
      1 &   1\\
      2 &   1+u\\
      3 &  1+u+2u^2+u^3+u^4\\
      4 &  1+u+2u^2+4u^3+4u^4+4u^5+4u^6+2u^7+u^8+u^{9}\\
      5 & 1+u+2u^2+4u^3+8u^4+11u^5+11u^6+14u^7+16u^{8}+\\
      & +14u^9+11u^{10}+11u^{11}+8u^{12}+4u^{13}+2u^{14}+u^{15}+u^{16}
    \end{array}
  \end{equation*}
  There is clearly a lot of structure in the polynomials $P_k(u)$
  which will be discussed in the next section.


\section{The structure of $P_k(u)$} 
Let us first write down the expression for $P_k(u)$ which came from
(\ref{ex}):
\begin{equation}\label{P}
  P_k(u) = \sum_{j=0}^k (-1)^{j}\Euler_{k-j}(u^{k+1}) (1-u^{k+1})^{j}
  \sum_{i=j}^k \binom{i}{j} u^{-i}.
\end{equation}
If we write $P_k(u) = \sum_{i=0}^{k^2} \alpha_i u^i$, we will define
the \emph{stretch} of $P_k(u)$ to be
\begin{equation*}
  PP_k(u) 
  = \alpha_0 + \sum_{i=0}^k \sum_{j=0}^{k-2} 
  \alpha_{1+i+(k+1)} u^{2+i+(k+1)j+j}.
\end{equation*}
What this does to $P_k(u)$ is to insert $0$ coefficients at every
$(k+1)$st term, starting after $\alpha_0$.  Thus, the stretched
polynomials corresponding to the values of $P_k(u)$ given in the array
above are:
\begin{equation*}
  \begin{array}{c|l}
    k & PP_k(u)  \\ \hline
    1 & 1\\
    2 & 1+u^2\\
    3 & 1+u^2+2u^3+u^4+u^6\\
    4 & 1+u^2+2u^3+4u^4+4u^5+4u^7+4u^8+2u^9+u^{10}+u^{12}\\
    5 & 1+u^2+2u^3+4u^4+8u^5+11u^6+11u^8+14u^9+16u^{10}+\\
      & +14u^{11}+11u^{12}+11u^{14}+8u^{15}+4u^{16}+2u^{17}+u^{18}+u^{20}
  \end{array}
\end{equation*}
Note that if $P_k(u)$ has degree $k^2$ then $PP_k(u)$ has degree
$k^2 + k$.  
\begin{theorem} \label{th2}
  For all $k \geq 1$, 
  \begin{equation*}
    P_{k+1}(u) = PP_k(u)\cdot (1 + u + u^2 + \dots + u^{k+1}).  
  \end{equation*}
\end{theorem}
\begin{proof}
  From (\ref{P}) and the definition of $PP_k(u)$, we can write
  \begin{equation*}
    PP_k(u) = \sum_{t=0}^k (-1)^t (1-u^{k+2})^t 
    \Euler_{k-t}(u^{k+2}) \sum_{s=t}^k \binom{s}{t} u^{-s}.
  \end{equation*}
  We want to show that 
  \begin{equation*}
    P_{k+1}(u) = PP_k(u) \cdot \frac{1-u^{k+2}}{1-u}.
  \end{equation*}
  Thus, it is enough to prove that $A = B$ where
  \begin{align*}
    A &= PP_{k+1}(u) (1-u^{k+2})  \\
      &= \sum_{t=0}^k (-1)^t (1-u^{k+2})^{t+1} 
    \Euler_{k-t}(u^{k+2})\sum_{s=t}^k\binom{s}{t} u^{-s}
    \shortintertext{and}
    B &= P_k(u) (1-u) \\
      &=  \sum_{t=0}^{k+1} (-1)^t (1-u^{k+2})^{t} 
    \Euler_{k+1-t}(u^{k+2})  \sum_{s=t}^{k+1} \binom{s}{t} u^{-s} (1-u).
  \end{align*}
  Observe that
  \begin{align*}
    B
    &= \sum_{t=0}^{k+1} (-1)^t (1-u^{k+2})^{t} \Euler_{k+1-t}(u^{k+2}
       \left( \sum_{s=t}^{k+1} \binom{s}{t} u^{-s} (1-u)\right)\\
    &=  \sum_{t=1}^{k+1} (-1)^t (1-u^{k+2})^{t} \Euler_{k+1-t}(u^{k+2})
       \left( \sum_{s=t}^{k+1} \binom{s}{t} u^{-s} (1-u)\right)\\
    & \quad + \Euler_{k+1}(u^{k+2}) (1-u^{-k-2}) (-u)\\
    &= \sum_{t=0}^{k} (-1)^t (1-u^{k+2})^{t+1} \Euler_{k-t}(u^{k+2})
       \left( \sum_{s=t+1}^{k}\binom{s}{t+1} u^{-s} (u-1)\right) \\
    & \quad - u (1-u^{-k-2}) \Euler_{k+1}\\
    &= \sum_{t=0}^{k} (-1)^t (1-u^{k+2})^{t+1} \Euler_{k-t}(u^{k+2}) 
       \left( \sum_{s=t+1}^{k} 
       \binom{s}{t} u^{-s} +u^{-t} - \binom{k+1}{t+1} u^{-k-1} \right) \\
    &\quad  - u (1-u^{-k-2}) \Euler_{k+1}\\
    &= \sum_{t=0}^{k} (-1)^t (1-u^{k+2})^{t+1} \Euler_{k-t}(u^{k+2}) 
       \left(
       \sum_{s=t}^{k} \binom{s}{t} u^{-s}-\binom{k+1}{t+1} u^{-k-1} 
       \right) \\
    & \quad - u (1-u^{-k-2}) \Euler_{k+1}.
  \end{align*}
  Hence, to prove that $A = B$, we only need to establish
  \begin{multline*}
    \sum_{t=0}^{k} (-1)^t (1-u^{k+2})^{t+1} \Euler_{k-t}(u^{k+2}) 
    \left( -\binom{k+1}{t+1} u^{-k-1} \right) 
    - u (1-u^{-k-2}) \Euler_{k+1}\\
    = u(1-u^{-k-2}) \Euler_{k+1}(u^{k+2}).
  \end{multline*}
  However, this would follow from
  \begin{equation}\label{euler-id}
    \sum_{t=0}^{k'} (x-1)^t \Euler_{k'-t}(x) \binom{k'}{t} = x \Euler_{k'}(x)
  \end{equation}
  by taking $k' = k+1$ and $x = u^{k+2}$. So it remains to prove
  (\ref{euler-id}).
  
  To do this we will use the standard generating function for
  $\Euler_n(w)$ from (\ref{euler_gf}):
  \begin{equation*}
    \sum_{n,m \geq 0} \euler{n}{m} w^m \frac{z^n}{n!} 
    = \sum_{n \geq 0} \Euler_n(w) \frac{z^n}{n!} 
    = \frac{1-w}{e^{(w-1)z} - w}.
  \end{equation*}
  Consider
  \begin{align*}
    F(x,z) &= \sum_{k > 0} x \Euler_k(x) \frac{z^k}{k!}
    \shortintertext{and}
    G(x,z) &= \sum_{k > 0} \sum_{t=0}^k (x-1)^t \Euler_{k-t} (x) \binom{k}{t}
    \frac{z^k}{k!}.
  \end{align*}
  It will suffice to show that $F = G$. Now
  \begin{align*}
    G(x,z) &= \sum_{k \geq 0} \sum_{t \geq 0} (x-1)^t \Euler_{k-t}(x)
    \binom{k}{t} \frac{z^k}{k!} - 1\\
    &= \sum_{t \geq 0} \sum_{k \geq t} (x-1)^t \Euler_{k-t}(x) 
    \binom{k}{t} \frac{z^k}{k!} - 1\\
    &= \sum_{t \geq 0} (x-1)^t \sum_{k' \geq 0}  \Euler_{k'}(x)
    \binom{k' + t}{t} \frac{z^{k' + t}}{(k' + t)!} - 1\\
    &= \sum_{t \geq 0} (x-1)^t z^t
    \sum_{k' \geq 0}\Euler_{k}(x) \frac{z^{k}}{k!} - 1\\
    &= e^{(x-1)z} \cdot \frac{1-x}{e^{(x-1)z} - x} - 1
    = \frac{x (1 - e^{(1-x)z})}{e^{(x-1)z} -x}.
  \end{align*}
  On the other hand,
  \begin{align*}
    F(x,z) 
    &= \sum_{k > 0} x \Euler_k(x) \frac{z^k}{k!}\\
    &= x \left( \frac{1-x}{e^{(x-1)z}-x} - 1 \right)
    =  \frac{x (1 - e^{(1-x)z})}{e^{(x-1)z} -x}
  \end{align*}
  as desired. This completes the proof of Theorem \ref{th2}. 
\end{proof}

\begin{theorem}
  The coefficients of $P_k(u)$ are symmetric and unimodal. 
\end{theorem}

\begin{proof}
  It follows from Theorem \ref{th2} that we can construct the
  coefficient sequence for $P_{k+1}(u)$ from that of $P_k(u)$ by the
  following rule (where we assume that all coefficients of $u^t$ in
  $P_k(u)$ are $0$ if $t < 0$ or $t > k^2$). Namely, suppose we write
  $P_k(u) = \sum_{i=0}^{k^2} \alpha_i u^i$ so that we have the
  coefficient sequence $A_k = (\alpha_0, \alpha_1, \dots,
  \alpha_{k^2})$. Now form the new sequence $B_k = (\beta_0, \beta_1,
  \dots \beta_{k^2 + k})$ by the rule
  \begin{equation*}
    \beta_i = \sum_{j=i-k}^i \alpha_j ,\quad 0 \leq i \leq k^2 + k.
  \end{equation*}
  Finally, starting with $\beta_0$, insert \emph{duplicate} values
  for the coefficients
  \begin{equation*}
  \beta_0,\beta_{k+1},\beta_{2(k+1)},\dots,\beta_{t(k+1)},
  \dots,\beta_{(k-1)(k+1)}\text{ and }\beta_{k(k+1)}.    
  \end{equation*}
  Thus, this will generate the sequence
  \begin{equation*}
  (\beta_0, \beta_0, \beta_1, \beta_2, \dots, \beta_k, \beta_{k+1},
  \beta_{k+1}, \beta_{k+2}, \dots, \beta_{k^2+k-1}, \beta_{k^2+k},
  \beta_{k^2+k}). 
  \end{equation*}
  This  new sequence will in fact just be the coefficient sequence
  $A_{k+1}$ for $P_{k+1}(u)$. For example, starting with $P_1(u) = 1 +
  u$, we have $A_1 = (1,1)$ and so $B_1 = (1, 2, 1)$. Now, inserting
  the duplicate values for $\beta_0 = 1$ and $\beta_2 = 1$, we get the
  coefficient sequence $A_2 = (1, {\bf 1}, 2, 1, {\bf 1})$ for $P_2(u)
  = 1 + u +2 u^2 + u^3 + u^4$.  Repeating this process for $A_2$, we
  sum blocks of length $3$ to get $B_2 = (1, 2, 4, 4, 4, 2,
  1)$. Inserting duplicates for entries at positions $0, 3$ and $6$
  gives us the new coefficient sequence $A_3 = (1, {\bf 1}, 2, 4, 4,
  {\bf 4}, 4, 2, 1, {\bf 1})$ of $P_3 = 1 + u + 2u^2 + 4u^3 +4u^4
  +4u^5 +4u^6 +2u^7 +u^8 +u^9$, etc. It is also clear from this
  procedure that if $A_k$ is symmetric and unimodal, then so is $B_k$,
  and consequently, so is $A_{k+1}$. This is what we claimed.
\end{proof}

\subsection{An Eulerian identity} 
Note that since $P_k(u)$ is symmetric and has degree $u^{k^2}$, we
have $P_k(u) = u^{k^2} P_k({\textstyle \frac{1}{u}})$.  Replacing
$P_k(u)$ by its expression in (\ref{P}), we obtain (with some
calculation) the interesting identity
\begin{equation*}
  \sum_{j = 0}^{a+b} (-1)^j \binom{a}{j} (1-x)^j \Euler_{a+b-j}(x) 
  = x \sum_{j = 0}^{a+b} \binom{b}{j}(1-x)^j \Euler_{a+b-j}(x)
  +\binom{b}{a+b}(1-x)^{a+b+1}
\end{equation*}
for \emph{all} integers $a$ and $b$ provided that $a + b \geq 0$.


\begin{thebibliography}{99}

\bibitem{jug0} J. Buhler, D. Eisenberg, R. Graham and C. Wright,
  Juggling drops and descents, {\it Amer. Math. Monthly} {\bf 101}
  (1994), 507--519.

\bibitem{jug} F. Chung, and R. L. Graham, Primitive juggling
  sequences, {\it Amer. Math. Monthly} {\bf 115}, March 2008,
  185--194.

\bibitem{euler} L. Euler, Methodus universalis series summandi
  ulterius promota, {\it Commentarii academiae scientiarum imperialis
    Petropolitanae} {\bf 8} (1736), 147--158. Reprinted in his {\it
    Pera Omnia}, series 1, volume 14, 124--137.

\bibitem{foata_han} D. Foata and G.-N. Han, $q$-series in Combinatorics;
  permutation statistics (Lecture Notes), preliminary edition, 2004.

\bibitem{foata_han_fix} D. Foata and G.-N. Han, Fix-Mahonian calculus III: A quadruple distribution,  {\it Monatsh. Math.}  {\bf{154}}  (2008),  no. 3, 177--197.

\bibitem{concrete} R. L. Graham, D. E. Knuth and O. Patashnik, {\it
  Concrete Mathematics}, Addison-Wesley, 1994.

\bibitem{knu} D.~E.~Knuth: {\it The Art of Computer Programming},
  Vol.~1, {\it Fundamental algorithms}. Addison-Wesley, Reading, 1969.

\bibitem{knuth} D.~E.~Knuth, {\it The Art of Computer Programming},
  Vol.~3,{\it Sorting and Searching}, Addison-Wesley, Reading, 2nd
  ed. 1998.

\bibitem{mac} P. A. MacMahon, {\it Combinatory Analysis}, 2 volumes,
  Cambridge University Press, London, 1915-1916. Reprinted by Chelsea,
  New York, 1960.

\bibitem{rodrigues} O. Rodrigues, Note sur les inversions, ou
  d\'erangements produits dans les permutations, {\it J. de Math.}
  {\bf{4}} (1839), 236--240.

\bibitem{wachs} J. Shareshian and M. L. Wachs, $q$-Eulerian
  polynomials: excedance number and major index.  {\it
    Electron. Res. Announc. Amer. Math. Soc.} {\bf 13} (2007), 33--45.

\bibitem{wachs2} J. Shareshian and M. L. Wachs, Eulerian
  quasisymmetric functions, preprint 2009.

\end{thebibliography}
\end{document}